\documentclass{amsart}
\usepackage{amssymb,latexsym,amsmath,amsthm,mathrsfs}
\usepackage{graphicx}

\begin{document}
\title{Euler and the pentagonal number theorem}
\thanks{The author
was supported by an NSERC USRA while this paper
was written.}
\author{Jordan Bell}
\email{jbell3@connect.carleton.ca}
\address{School of Mathematics and Statistics, Carleton
University, 1125 Colonel By Drive, Ottawa, Ontario, K1S 5B6, Canada.}
\date{\today}

\begin{abstract}
In this paper we give the history of Leonhard Euler's work
on the pentagonal number theorem, and his applications of the pentagonal
number theorem to the divisor function, partition function and divergent series. We have
attempted to give an exhaustive review of all of Euler's correspondence
and publications about the pentagonal number theorem and his applications
of it.

Comprehensus: In hoc dissertatione damus historiam operis Leonhardi Euleri
super theorma numerorum pentagonalium, et eius usus theoremae
numerorum pentagonalium ad functioni divisori, functioni partitione et
seriebus divergentibus. Conati sumus dedisse recensum plenum omnium
commerciorum epistolicarum et editionum Euleri circa theorema
numerorum pentagonalium et eius applicationum ipsius.
\end{abstract}
\maketitle

\section{Introduction}
\label{section:introduction}
The pentagonal numbers are those numbers of the form $\frac{n(3n-1)}{2}$
for $n$ a positive integer.
They represent the number of distinct points which may be arranged
to form superimposed regular pentagons with the same number of equally spaced
points on the sides of each respective pentagon. The generalized
pentagonal numbers are those numbers of the form $\frac{n(3n \pm 1)}{2}$ for
$n$ non-negative,
i.e. the pentagonal numbers for $n$ an integer.
However, in this paper unless we indicate otherwise, by pentagonal number
we will mean a generalized pentagonal number.

In this paper we will consider Leonhard Euler's work on the pentagonal number
theorem and his applications of it to recurrence relations for 
the divisor function and the partition function, and to divergent series.
We have attempted to give an exhaustive summary of Euler's
correspondence and works that discuss the pentagonal number theorem.
The pentagonal number theorem is the
formal identity:
\begin{equation}
\label{equation:maintheorem}
\prod_{m=1}^\infty (1-x^m) = \sum_{n=-\infty}^\infty (-1)^n x^{\frac{n(3n-1)}{2}},
\end{equation}
and it is called the pentagonal number theorem because the exponents
in the formal power series on the right-hand side of the equation
are the pentagonal numbers. Formally,
the set of all power series over the integers form a commutative
$\mathbf{Z}$-algebra
closed under logarithmic differentiation and integration.
Analytically, the series on the right-hand side of \eqref{equation:maintheorem}
converges absolutely for $|x|<1$; in fact, for the function
$\phi$ defined by
$\phi(x)=\prod_{m=1}^\infty (1-x^m)$ for $|x|<1$, Konrad Knopp \cite{knopp}
shows that $\phi$ does not have an analytic continuation beyond the unit
circle.

We can observe in fact that the identity \eqref{equation:maintheorem} follows
directly from the Jacobi triple product
identity $\prod_{m=1}^\infty (1-q^{2m})(1+q^{2m-1}z^2)(1+q^{2m-1}z^{-2})=
\sum_{n=-\infty}^\infty q^{n^2}z^{2n}$, for $q=x^{\frac{3}{2}}$ and
$z^2=-x^{1/2}$. This gives
$\prod_{m=1}^\infty (1-x^{3m})(1-x^{3m-1})(1-x^{3m-2})=\sum_{n=-\infty}^\infty
(-1)^n x^{3n(n+1)}{2}$, that is 
\eqref{equation:maintheorem}.
 We further discuss the Jacobi triple product identity
 in Section~\ref{section:conclusion}, and also 
 $q$-series and theta functions. 
Fabian Franklin \cite{franklin} gives a combinatorial proof of
the pentagonal number theorem
by proving an equivalent result about partitions.
George E. Andrews \cite{andrews} generalizes Euler's proof of the pentagonal
number theorem to prove that $1-\sum_{m=1}^\infty
(1-zq)(1-zq^2)\cdots(1-zq^{m-1})z^{m+1}q^m=1+\sum_{n=1}^\infty
(-1)^n\big(z^{3n-1}q^{\frac{n(3n-1)}{2}}+z^{3n}q^{\frac{n(3n+1)}{2}}\big)$;
analytically, these series converge absolutely for $|q|<1$ and $|z|<|q|^{-1}$.

Following Euler's notation, in this paper we express by $\int n$ the divisor
function, the sum of
all the divisors of $n$ including itself, and by
$n^{(m)}$ the partition function, the number 
of ways of writing $n$ as a sum of positive integers less than or 
equal to $m$, disregarding order.

\section{Background}
\label{section:background}
The first time the pentagonal number theorem is mentioned in 
Euler's correspondence is in a letter from Daniel I Bernoulli
to Euler on January 28, 1741, Letter XX in the Daniel I Bernoulli-Euler
correspondence in P.-H. Fu{\ss}'s collection
{\em Correspondance math\'ematique et physique de quelques
              c\'el\`ebres g\'eom\`etres du {XVIII}\`eme si\`ecle}
\cite{fuss1843}, OO140 in the annotated index of Euler's correspondence
in the {\em Opera omnia} \cite{commercium1}.
In this letter Bernoulli discusses a number of problems
that Euler apparently had posed to him in earlier letters (the letters from Euler in the Daniel Bernoulli I-Euler
correspondence are not extant); in particular, he mentions the problem
of finding all the partitions of an integer. Concerning the pentagonal
number theorem, 
Bernoulli writes: ``The other problem, to transform the expression
$\big(1-\frac{1}{n}\big)\big(1-\frac{1}{n^2}\big)\big(1-\frac{1}{n^3}\big)$ into the
series
$1-\frac{1}{n}-
\frac{1}{n^2}+\frac{1}{n^5}+\frac{1}{n^7}-\frac{1}{n^{12}}-\frac{1}{n^{15}}+$ etc.
follows easily by induction,
if one multiplied many factors from the
given expression. The remaining of the series, in which prime numbers are
seen, I do not see. This can be shown in a most pleasant investigation,
together with tranquil pastime and the endurance of pertinacious labor, all
three of which I lack''. 
Thus Euler probably had mentioned the problem of expanding this infinite
product into an infinite series in his last letter, on September 15, 1740.

The first paper in which Euler mentions the pentagonal number theorem
is his  ``Observationes analyticae variae de combinationibus'',
presented to the St. Petersburg Academy on April 6, 1741 and published
in 1751 in the {\em Commentarii academiae scientiarum imperialis Petropolitanae}
\cite{E158}, E158 in the Enestr\"om index \cite{enestrom13}. In this paper Euler
introduces the generating function
\begin{equation}
\label{equation:partition}
\prod_{x=1}^\infty \frac{1}{1-n^x}
=1+n+2n^2+3n^3+5n^4+7n^5+11n^6+15n^7+22n^8+\textrm{etc.}
\end{equation}
for the (unrestricted) partition function, i.e. $n^{(\infty)}$.
 In \S 36 of
this paper, Euler says, ``Here at the end of this dissertation
a noteworthy observation should be made, which however I have not been able to
demonstrate with geometric rigor. I have observed
namely for this infinite product:
\[
(1-n)(1-n^2)(1-n^3)(1-n^4)(1-n^5)\textrm{ etc.},
\]
if expanded by multiplication, to produce this series:
\[
1-n-n^2+n^5+n^7-n^{12}-n^{15}+n^{22}+n^{26}-n^{35}-n^{40}+n^{51}+\textrm{etc.},
\]
where each of these occurs as a power of $n$, of which the exponents
are contained in the form $\frac{3xx \pm x}{2}$. And if
$x$ is an odd number, the powers of $n$, which are $n^{\frac{3xx \pm x}{2}}$,
will have the coefficient $-1$; and if $x$ is an even number, then the 
powers $n^{\frac{3xx \pm x}{2}}$ will have the coefficient $+1$''.
In \S 37, Euler notes as well that the product of the power series
on the right-hand side of \eqref{equation:partition} and
the above power series
is unity, since they
are the series expansions of reciprocal infinite products.

The next time the expansion of the infinite product
$\prod_{m=1}^\infty (1-x^m)$
into a series comes up in Euler's correspondence is in a letter from
Euler to Niklaus I Bernoulli 
on September 1, 1742, Letter 3 in the Euler-Niklaus I Bernoulli
correspondence in the {\em Opera omnia} \cite{commercium2}, OO236.
However, Euler does discuss the problem of expressing the series $s-\frac{s^3}{6}+
\frac{s^5}{120}-\textrm{etc.}$ as the infinite 
product $s(1-\frac{ss}{\pi \pi})(1-\frac{ss}{4\pi
\pi})(1-\frac{ss}{9\pi \pi})$ etc.
before the pentagonal number theorem in this letter,
and this is also discussed in the previous two letters in the
Euler-Niklaus I Bernoulli correspondence. Euler also
notes in this letter that the coefficients of the terms in the series
\[
1+1n+2n^2+3n^3+5n^4+7n^5+11n^6+15n^7+22n^8+30n^9+42n^{10}+56n^{11}+\textrm{etc.}
\]
give the number of different ways in which the exponent of the term
can be made by addition, i.e. that it is the generating function for the
(unrestricted) partition function. Euler then writes: ``This series moreover
arises from division, if unity were divided by
$(1-n)(1-n^2)(1-n^3)(1-n^4)(1-n^5)$ etc., which product if expanded
gives this expression
\[
1-n-n^2+n^5+n^7-n^{12}-n^{15}+n^{22}+n^{26}-n^{35}-\textrm{etc.}
\]
where the precise way in which the exponents proceed I have not been able 
to penetrate, although by induction I have concluded for no other
exponents to occur, unless they are contained in the formula
$(3xx \pm x)/2$; and this is such that the powers of $n$ have the $+$ sign
if the exponents arise with an even number substituted for $x$''.

Niklaus I Bernoulli replies to Euler on October 24, 1742, Letter 4 in the Euler-Niklaus
I Bernoulli correspondence in the {\em Opera omnia} \cite{commercium2},
OO237.
Bernoulli discusses the pentagonal number theorem after
mentioning the generating function $\prod_{m=1}^\infty \frac{1}{1-x^m}$ for the partition function.
Bernoulli here observes, ``In the expansion of the series
\[
n^0-n^1-n^2+n^5+n^7-n^{12}-n^{15}+n^{22}+n^{26}-n^{35}-\textrm{etc.}
\]
which you have found to be equal to the product $(1-n)(1-nn)(1-n^3)$ etc.,
the differences of the exponents progress as $1,1,3,2,5,3,7,4,9,5$, etc.,
which numbers taken alternately are from the series $1,3,5,7,9$, etc.
and from the series $1,2,3,4,5$, etc., which properties will perhaps
be able to be demonstrated from the nature of this thing not only
through induction; but into this matter it is not now free to inquire''.

Euler next writes to Niklaus I Bernoulli on November 10, 1742,
Letter 5 in the Euler-Niklaus I Bernoulli correspondence 
in the {\em Opera omnia} \cite{commercium2}, OO238, the last
letter in their correspondence that deals with the pentagonal
number theorem.
Euler writes: ``This expression
\[
(1-n)(1-n^2)(1-n^3)(1-n^4)\textrm{ etc.}
\]
by expansion shall give the series
\[
1-n-n^2+n^5+n^7-\textrm{etc.}
\]
in which no other exponents occur unless they are contained
in
\[
\frac{3xx \pm x}{2},
\]
which I have for my part concluded with legitimate induction,
even if I have not been able to find a demonstration in any manner,
however, although I have not devoted enough time to this. 
I however have found this expression
\[
(1-n)(1-n^2)(1-n^3)(1-n^4)\textrm{ etc.}
\]
to be able indeed to be transformed into this series
\[
1-\frac{n}{1-n}+\frac{n^3}{(1-n)(1-n^2)}-\frac{n^6}{(1-n)(1-n^2)(1-n^3)}+
\textrm{etc.}
\]
of which the value is equal precisely to the sum of the series
\[
1-n^1-n^2+n^5+n^7-n^{12}-n^{15}+\textrm{etc.''}
\]
In this letter Euler also gives more general results
on the infinite products $(1+mz)(1+m^2z)(1+mz^3)(1+mz^4)$ etc.
and
$\frac{1}{(1-mz)(1-m^2z)(1-m^3z)(1-m^4z)\textrm{ etc.}}$. Euler shows
that
\begin{eqnarray*}
(1+mz)(1+m^2z)(1+m^3z)(1+m^4z)\textrm{ etc.}\\
=1+\frac{mz}{1-m}+\frac{m^3z^2}{(1-m)(1-m^2)}+
\frac{m^6z^3}{(1-m)(1-m^2)(1-m^3)}+\textrm{ etc.}
\end{eqnarray*}
and
\begin{eqnarray*}
\frac{1}{(1-mz)(1-m^2z)(1-m^3z)\textrm{etc.}}\\
=1+\frac{mz}{1-m}+\frac{m^2z^2}{(1-m)(1-m^2)}+
\frac{m^3z^3}{(1-m)(1-m^2)(1-m^3)}+\textrm{ etc.}
\end{eqnarray*}

The next time Euler discusses the pentagonal number theorem in his
correspondence is in a letter to Christian Goldbach on October 15, 1743,
Letter 74 in the Euler-Goldbach correspondence \cite{eulergoldbach},
OO788. Euler states here,
``If these factors $(1-n)(1-n^2)(1-n^3)(1-n^4)(1-n^5)$ etc.
are multiplied out onto infinity, the following series 
\[
1-n^1-n^2+n^5+n^7-n^{12}-n^{12}+n^{22}+n^{26}-n^{35}-n^{40}+n^{51}+n^{57}-
\textrm{etc.}
\]
is produced, 
from which it is easily shown by induction that all of the terms
are contained in the form $n^\frac{3xx \pm x}{2}$, and that they have the prefixed
sign $+$ when $x$ is an even number, and the sign $-$ when $x$ is odd.
I have however not yet found a method by which I could prove the identity
of these two expressions. The Hr. Prof. Niklaus Bernoulli has also been able to prove
nothing beyond induction''.  

Goldbach replies to Euler's problem in a letter dated December 1743,
Letter 75 in the Euler-Goldbach collection \cite{eulergoldbach},
OO789. Goldbach does not give any explicit
ideas for how to tackle proving the pentagonal number theorem, and instead poses
a new related problem.
He says, ``About the 
series $(1-n)(1-nn)(1-n^3)$ etc., another problem has occurred to me:
Given an infinite series of terms $A$, with an order given for the varying
signs
$+$ and $-$ of the progression, to find a series $B$ of such a nature that
in the product $AB$ the signs $+$ and $-$ succeed in the same order
that they succeeded in $A$.
This problem can easily be solved in the case $A=(1-n)(1-nn)(1-n^3)$ etc., 
although as it has already been noted the signs $+$ and $-$ alternate
in an unusual manner, for if I set
$B=(1-n^\frac{1}{2})(1-n^\frac{3}{2})(1-n^\frac{5}{2})$ etc.,
then $A$ multiplied by $B$ becomes a new series, which contains
the same  variation of
the signs''.

Euler writes back to Goldbach about this on January 21, 1744,
Letter 76 in the Euler-Goldbach correspondence \cite{eulergoldbach},
OO790,
in which he says, ``Reflection about the expression
$(1-n)(1-n^2)(1-n^3)$ etc., in view of the factors
$(1-n^\frac{1}{2})(1-n^\frac{3}{2})(1-n^\frac{5}{2})$, which,
having been composed, if expanded, yields an equal alternation
of the signs $+$ and $-$, could perhaps be advantageous in other research;
alone in the series, which I derived from this, I have 
not been able to make any use of it''. 

On April 5, 1746, Euler writes another letter to Goldbach
about the pentagonal number theorem,
Letter 102 in the Euler-Goldbach correspondence
\cite{eulergoldbach},
OO816.
In this letter he considers several problems, in particular
the expansion of the infinite products
$(1+a)(1+a^2)(1+a^4)(1+a^8)\textrm{ etc.}$
and $(1-a)(1-a^2)(1-a^4)(1-a^8)\textrm{ etc.}$ into series. Euler
states as a theorem that ``If it were
$s=(1-na)(1-n^2a)(1-n^3a)(1-n^4a)(1-n^5a)$ etc. onto infinity, it will be
\begin{align*}
s=1-\frac{na}{1-n}+\frac{n^3a^2}{(1-n)(1-n^2)}-
\frac{n^6a^3}{(1-n)(1-n^2)(1-n^3)}\\
+\frac{n^{10}a^4}{(1-n)(1-n^2)(1-n^3)(1-n^4)}-\textrm{etc.}
\end{align*}
and
\begin{align*}
\frac{1}{s}=1+\frac{na}{1-n}+\frac{n^2a^2}{(1-n)(1-n^2)}+
\frac{n^3a^3}{(1-n)(1-n^2)(1-n^3)}\\
+\frac{n^4a^4}{(1-n)(1-n^2)(1-n^3)(1-n^4)}+\textrm{etc.}
\end{align*}
I believe I have also written that if one multiplies
this product onto infinity
\[
(1-a)(1-a^2)(1-a^3)(1-a^4)(1-a^5)\textrm{ etc.}
\]
this series is produced
\[
1-a-a^2+a^5+a^7-a^{12}-a^{15}+a^{22}+a^{26}-a^{35}-a^{40}+\textrm{etc.},
\]
where the order of the exponents is very peculiar, and also by
induction it may be determined that all are contained in the form
$\frac{3xx \pm x}{2}$, though I have not yet been able to expose
the rule which has been observed from the nature of this matter''.

Later, in a letter to Jean le Rond d'Alembert on December 30, 1747, 
 Letter 11 in the Euler-d'Alembert correspondence in the {\em Opera omnia}
\cite{commercium5}, OO23,
Euler writes that he has learned from de Maupertuis (the President
of the Berlin Academy) that d'Alembert wants to leave
his mathematical research for some time to regain his health. 
Euler goes on to say, ``If in your spare time you should wish to do some
research which does not require much effort, I will take liberty
to propose the expression $(1-x)(1-x^2)(1-x^3)(1-x^4)(1-x^5)(1-x^6)$ etc.,
which upon expansion by multiplication gives the series
\[
1-x^1-x^2+x^5+x^7-x^{12}-x^{15}+x^{22}+x^{26}-x^{35}-x^{40}+x^{51}+x^{57}-
x^{70}-x^{77}+\textrm{etc.}
\]
which would seem very remarkable to me because of the law which we easily
discover within it, but I do not see how this law may be deduced
without induction of the proposed expression''. 
As a postscript to the letter, Euler writes,
``If we put $s=(1-x)(1-x^2)(1-x^3)(1-x^4)(1-x^5)(1-x^6)$ etc., I can show that
\begin{align*}
s=1-\frac{x}{1-x}+\frac{x^3}{(1-x)(1-x^2)}-\frac{x^6}{(1-x)(1-x^2)(1-x^3)}\\
+\frac{x^{10}}{(1-x)(1-x^2)(1-x^3)(1-x^4)}-\textrm{etc.''}
\end{align*}

d'Alembert replies to Euler 
in a letter on January 20, 1748,
Letter 12 in the Euler-d'Alembert correspondence in the {\em Opera
omnia} \cite{commercium5}, OO24.
d'Alembert says, ``regarding the series of which you have spoken, it is
very peculiar, and I have reflected on it for a while, but I only
see induction to show it. At the end, no one is deeper and better versed
on these matters than you''.

In \S 323 (in Chapter XVI, ``De partitione numerorum'') of the {\em Introductio in analysin
infinitorum} \cite{E101}, E101, published in 1748,
Euler notes
in his discussion on the generating
function $\prod_{m=1}^\infty \frac{1}{1-x^m}$
of the partition function
that, ``in particular it should be observed for there to be a ladder
relationship in the denominator, for if
indeed the factors of the denominator are multiplied successively
into each other, it will advance
\[
1-x-x^2+x^5+x^7-x^{12}-x^{15}+x^{22}+x^{26}-x^{35}-x^{40}+x^{51}+\textrm{etc.},
\]
in which series if it were considered attentively, for no other powers of $x$
to be discovered, unless the exponents of them are contained in the formula
$\frac{3nn \pm n}{2}$; and too, if $n$ is an odd number,
the powers will be negative; positive if indeed $n$ is an even number''.

\section{Euler's proof of the pentagonal number theorem}
\label{section:pentagonal}
In a letter from Goldbach to Euler on April 15, 1747, Letter 114
in the Euler-Goldbach correspondence \cite{eulergoldbach},
OO828,
Goldbach responds to Euler's previous letter that had given
a recurrence relation
for the divisor function (cf. Section~\ref{section:divisor}).
In that letter, Euler had remarked that his proof assumes the pentagonal
number theorem, which he had not been able to rigorously prove.
Goldbach declares that, ``The observation which you have communicated to me seems
to me already through the given induction proved to the extent
that one can believe in its truth one-hundred to one. 
Moreover, it has already been noted earlier that
$A\dots (1-x)(1-xx)(1-x^3)(1-x^4)\textrm{ etc.}=
B\dots 1-x-xx+x^5+x^7-x^{12}-x^{15}+$ etc.
and I remember that from this I came to the simple conclusion that
when the powers of $x$ in $B$ are doubled, and with
\[
C=1-xx-x^4+x^{10}+x^{14}-x^{24}-x^{30}+\textrm{etc.}
\]
then it must be
\[
\frac{C}{B}=(1-x)(1-x^3)(1-x^5)(1-x^7)\textrm{ etc.''}
\]

Euler replies to Goldbach on May 6, 1747, Letter 115 in the Euler-Goldbach
correspondence \cite{eulergoldbach}, OO829, saying:
``The observation which was made about the identity
\begin{align*}
A\dots (1-x)(1-x^2)(1-x^3)(1-x^4)\textrm{ etc.}=\\
B\dots 1-x-x^2+x^5+x^7-\textrm{etc.}
\end{align*}
that if
\[
C=1-x^2-x^4+x^{10}+x^{14}-x^{24}-x^{30}+\textrm{etc.},
\]
then $\frac{C}{B}=(1-x)(1-x^3)(1-x^5)$ etc., I still remember well.
I have however neither from this nor from other considerations 
been able to display the identity between the formulas $A$ and $B$ properly;
for the fact that $A=B$ and that the exponents of $x$ in $B$ continue
only according to the series, 
I have also only been able to conclude by induction, which I
however have continued so far, that I consider the matter completely
true; I would be very enthusiastic to see a direct proof of this matter,
which would certainly lead to the discovery of many other beautiful
properties
of numbers; hitherto all of my pains have been for nothing''.

At last in a letter from Euler to Goldbach on June 9, 1750, Letter 144
in the Euler-Goldbach correspondence \cite{eulergoldbach},
OO858,
Euler
gives a proof of the pentagonal number theorem. Euler recalls
his discovery of a recurrence relation for the divisor function, and
that it assumes the pentagonal number theorem, which he had not been able
to prove. He then says: ``Since that time however I have also found a demonstration
of this theorem, which depends on this lemma:
\begin{align*}
(1-\alpha)(1-\beta)(1-\gamma)(1-\delta)\textrm{ etc.}=\\
1-\alpha-\beta(1-\alpha)-\gamma(1-\alpha)(1-\beta)-\delta(1-\alpha)(1-\beta)(1-\gamma)-\textrm{etc.},
\end{align*}
whose demonstration is immediate.

Therefore according to this lemma
\begin{align*}
(1-x)(1-x^2)(1-x^3)(1-x^4)(1-x^5)\textrm{ etc.}=s
=\\
1-x-x^2(1-x)-x^3(1-x)(1-x^2)-x^4(1-x)(1-x^2)(1-x^3)-\textrm{etc.}
\end{align*}
Were it put $s=1-x-Axx$, it will be
\[
A=1-x+x(1-x)(1-x^2)+x^2(1-x)(1-x^2)(1-x^3)+\textrm{etc.}
\]
The factor $1-x$ is expanded everywhere, and it will be
\begin{align*}
A&=1&-x&-x^2(1-xx)&-x^3(1-xx)(1-x^3)&-\textrm{etc.}\\
&&+x(1-xx)&+x^2(1-xx)(1-x^3)&+x^3(1-xx)(1-x^3)(1-x^4)&+\textrm{etc.}
\end{align*}
and then it can be made
\[
A=1-x^3-x^5(1-xx)-x^7(1-xx)(1-x^3)-\textrm{etc.}
\]
Were it $A=1-x^3-Bx^5$, it will be
\[
B=1-xx+x^2(1-xx)(1-x^3)+x^4(1-xx)(1-x^3)(1-x^4)+\textrm{etc.}
\]
The factor $1-xx$ is expanded:
\begin{align*}
B&=1&-xx&-x^4(1-x^3)&-x^6(1-x^3)(1-x^4)&-\textrm{etc.}\\
&&+xx(1-x^3)&+x^4(1-x^3)(1-x^4)&+x^6(1-x^3)(1-x^4)(1-x^5)&+\textrm{etc.},
\end{align*}
and then it can be made
\[
B=1-x^5-x^3(1-x^3)-x^{11}(1-x^3)(1-x^4)-\textrm{etc.}
\]
Were it $B=1-x^5-Cx^8$, it will be
\[
C=1-x^3+x^3(1-x^3)(1-x^4)+x^6(1-x^3)(1-x^4)(1-x^5)+\textrm{etc.}
\]
The factor $1-x^3$ is expanded:
\begin{align*}
C&=1&-x^3&-x^6(1-x^4)&-x^9(1-x^4)(1-x^5)&-\textrm{etc.}\\
&&+x^3(1-x^4)&+x^6(1-x^4)(1-x^5)&+x^9(1-x^4)(1-x^5)(1-x^6)&+\textrm{etc.},
\end{align*}
therefore
\[
C=1-x^7-x^{11}(1-x^4)-x^{15}(1-x^4)(1-x^5)-\textrm{etc.}
\]
Were it $C=1-x^7-Dx^{11}$ etc. If one continues in the same manner, then
\[
D=1-x^9-Ex^{14}, \quad E=1-x^{11}-Fx^{17}, \quad \textrm{etc.}
\]
And in consequence:
\[
\left.\begin{array}{l}
s=1-x-Ax^2\\
A=1-x^3-Bx^5\\
B=1-x^5-Cx^8\\
C=1-x^7-Dx^{11}\\
D=1-x^9-Ex^{14}\\
\textrm{etc.}
\end{array}\right\}
\textrm{ or }
\left\{\begin{array}{l}
s=1-x-Ax^2\\
Ax^2=x^2(1-x^3)-Bx^7\\
Bx^7=x^7(1-x^5)-Cx^{15}\\
Cx^{15}=x^{15}(1-x^7)-Dx^{26}\\
Dx^{26}=x^{26}(1-x^9)-Ex^{40}\\
\textrm{etc.}
\end{array}\right.
\]
from which it doubtlessly follows
\[
s=1-x-x^2(1-x^3)+x^7(1-x^5)-x^{15}(1-x^7)+x^{26}(1-x^9)-\textrm{etc.}
\]
or
\[
s=1-x-x^2+x^5+x^7-x^{12}-x^{15}+x^{22}+x^{26}-x^{35}-\textrm{etc.''}
\]

In his ``Demonstratio theorematis circa ordinem in summis divisorum
observatum'', published in the {\em Novi commentarii academiae scientiarum imperialis
Petropolitanae} in 1760 \cite{E244}, E244, Euler 
gives the above proof of the pentagonal number theorem from
 his letter to Goldbach in more detail.
Euler first recalls that some time ago he had discovered a recurrence relation
for the divisor function, for which the differences in the arguments
are the pentagonal numbers, but that his proof was not rigorous enough
to him, since it relied on the pentagonal number theorem. But now
he 
declares that he has finally found 
a demonstration of this, and thus ``Not any doubt may remain of
this property, the demonstration of the truth of which rests
on each of these propositions, which I will now set forth and demonstrate''.
In Proposition~I, Euler notes that for
$s=(1+\alpha)(1+\beta)(1+\gamma)(1+\delta)(1+\varepsilon)(1+\zeta)(1+\eta)$
etc., then
\begin{align*}
s=(1+\alpha)+\beta(1+\alpha)+\gamma(1+\alpha)(1+\beta)+\delta(1+\alpha)(1+\beta)(1+\gamma)\\
+\varepsilon(1+\alpha)(1+\beta)(1+\gamma)(1+\delta)+
\zeta(1+\alpha)(1+\beta)(1+\gamma)(1+\delta)(1+\varepsilon)+\textrm{etc.}
\end{align*}
Taking $\alpha=-x$, $\beta=-xx$, $\gamma=-x^3$, $\delta=-x^4$,
$\varepsilon=-x^5$, etc., Euler then proves Proposition~II, that
if $s=(1-x)(1-xx)(1-x^3)(1-x^4)(1-x^5)(1-x^6)$ etc., then
\[
s=1-x-xx(1-x)-x^3(1-x)(1-x^2)-x^4(1-x)(1-x^2)(1-x^3)-\textrm{etc.}
\]
In Proposition~III, Euler then proves the pentagonal number
theorem, giving precisely the same proof which he gave in his above
letter
to Goldbach.
Euler then uses this to finally give a rigorous proof of his recurrence relation
for the divisor function, which we discuss in Section~\ref{section:divisor}.

Euler gives two proofs of the pentagonal number theorem
in ``Evolutio producti infiniti $(1-x)(1-xx)(1-x^3)(1-x^4)(1-x^5)(1-x^6)$ etc.
in seriem simplicem'', delivered to the St. Petersburg Academy on August 14, 1775 and published in 1783 in the {\em Acta academiae scientiarum imperialis Petropolitinae} \cite{E541}, E541. 
The first proof in \S \S 1--10 is the same as his earlier proof
except that instead of taking $s=1-x-Axx$, Euler takes $s=1-x-A$, 
instead of taking $s=1-x^3-Bx^5$, Euler takes $A=xx-x^5-B$, instead
of taking $B=1-x^5-Cx^8$, Euler takes $B=x^7-x^{12}-C$, etc.
In the second proof in \S \S 11--18, having taken $s=(1-x)(1-x^2)(1-x^3)$ etc., Euler states,
``Of course then
\[
s=1-x-xx(1-x)-x^3(1-x)(1-xx)-x^4(1-x)(1-x^2)(1-x^3)-\textrm{etc.},
\]
for which by expanding the second member $-xx(1-x)$, it will be
\[
s=1-x-xx+x^3-x^3(1-x)(1-xx)-x^4(1-x)(1-xx)(1-x^3)-\textrm{etc.},
\]
and it is set $s=1-x-xx+A$ so that
\[
A=x^3-x^3(1-x)(1-xx)-x^4(1-x)(1-xx)(1-x^3)-\textrm{etc.},
\] 
for which, by expanding each term by the factor $1-x$, it is broken into
two parts such that it appears as
\begin{align*}
A=x^3-x^3(1-xx)-x^4(1-xx)(1-x^3)-x^5(1-x^2)(1-x^3)(1-x^4)\ldots\\
+x^4(1-xx)+x^5(1-xx)(1-x^3)+x^6(1-x^2)(1-x^3)(1-x^4)\ldots
\end{align*}
Here, contracting the pairs of terms with the same powers of $x$ like before will
produce this
\begin{displaymath}
A=+x^5+x^7(1-xx)+x^9(1-xx)(1-x^3)+x^{11}(1-x^2)(1-x^3)(1-x^4)\textrm{ etc.''}
\end{displaymath}
Euler then expands again the second term, and takes $A=x^5+x^7-B$.
He expands each term in $B$ then by the factor $1-xx$, and expresses
this as for $A$ in a series of terms with negative signs and a series
of terms with positive signs, which he then combines to get
$B=x^{12}+x^{15}(1-x^3)+x^{18}(1-x^3)(1-x^4)+x^{21}(1-x^3)(1-x^4)(1-x^5)+
\textrm{etc.}$. Euler continues this up to $E=x^{51}+x^{57}(1-x^6)+
x^{63}(1-x^6)(1-x^7)+x^{69}(1-x^6)(1-x^7)(1-x^8)+\textrm{etc.}$,
and notes that, ``Here indeed the order of the exponents is easily
perceived. For with the values of the letters $A,B,C,D$,
at first the first terms were simply $x^3,x^9,x^{18},x^{30},x^{45}$,
where the exponents are clearly triples of triangular numbers,
from which in general the exponent for the number $n$ will be
$\frac{3nn+3n}{2}$. However, the differences between two
successive powers of $x$ is the same as $n$, because of which the number
$n$ will be subtracted twice from this formula, from which the exponents
follow as $\frac{3nn+n}{2}$ and $\frac{3nn-n}{2}$''.

We now shall give a formal proof by induction of the pentagonal number theorem
that is essentially Euler's above proof. It will be convenient to write
$\omega_n=\frac{3n(n-1)}{2}$.
Let $P=\prod_{m=1}^\infty
(1-x^m)=(1-x)(1-x^2)(1-x^3)(1-x^4)(1-x^5)\cdots$.
Thus
\begin{eqnarray*}
P&=&1-x-(1-x)x^2-(1-x)(1-x^2)x^3-(1-x)(1-x^2)(1-x^3)x^4\\
&&-(1-x)(1-x^2)(1-x^3)(1-x^4)x^5-\ldots
\end{eqnarray*}

Expanding out the factor $1-x$ in each term
 of $P$ yields
\begin{eqnarray*}
P&=&1-x-x^2+x^3-(1-x^2)x^3+(1-x^2)x^4-(1-x^2)(1-x^3)x^4\\
&&+(1-x^2)(1-x^3)x^5-(1-x^2)(1-x^3)(1-x^4)x^5\\
&&+(1-x^2)(1-x^3)(1-x^4)x^6-+\ldots
\end{eqnarray*}

Adding the coefficients of each power of $x$ gives
\begin{eqnarray*}
P&=&1-x-x^2+x^5+(1-x^2)x^7+(1-x^2)(1-x^3)x^9\\
&&+(1-x^2)(1-x^3)(1-x^4)x^6-+\ldots
\end{eqnarray*}

Since all the following terms have degree $\geq 6$,
it is clear that the first three terms of $P$ are $1-x-x^2$.
We now make the induction assumption that up to some $k$,
the first $2k+1$ terms of $P$ satisfy
\begin{eqnarray*}
P&=&\sum_{n=-k+1}^k (-1)^n x^{\omega_n}\\
&&+(-1)^k x^{\omega_k+k}\bigg((1-x^k)+(1-x^k)(1-x^{k+1})x^k\\
&&+(1-x^k)(1-x^{k+1})(1-x^{k+2})x^{2k}+\ldots\bigg)
\end{eqnarray*}
Clearly this holding for arbitrary $k$ is equivalent to the pentagonal number
theorem.

Now expanding out the factor $1-x^k$ in each term we obtain
\begin{eqnarray*}
P&=&\sum_{n=-k+1}^k (-1)^n x^{\omega_n}\\
&&+(-1)^kx^{\omega_k+k}\Big(1-x^k+(1-x^{k+1})x^k-(1-x^{k+1})x^{2k}\\
&&+(1-x^{k+1})(1-x^{k+2})x^{2k}-(1-x^{k+1})(1-x^{k+2})x^{3k}+-\ldots\Big)
\end{eqnarray*}

Adding the coefficients of each power then gives
\begin{eqnarray*}
P&=&\sum_{n=-k+1}^k (-1)^n x^{\omega_n}\\
&&+(-1)^kx^{\omega_k+k}\Big(1-x^{2k+1}+(1-x^{k+1})x^{3k+2}-(1-x^{k+1})(1-x^{k+2})x^{3k}+-\ldots\Big)
\end{eqnarray*}

Therefore we have
\begin{eqnarray*}
P&=&\sum_{n=-k+1}^k (-1)^n
x^{\omega_n}+(-1)^kx^{\omega_k+k}+(-1)^{k+1}x^{\omega_k+3k+1}\\
&&+(-1)^{k+1}x^{\omega_k+4k+2}\Big(1+(1-x^{k+1})x^{k+1}+(1-x^{k+1})(1-x^{k+2})x^{2k+2}+\ldots\Big)
\end{eqnarray*}

But $\omega_k+k=\frac{k(3k-1)+2k}{2}=\frac{k(3k+1)}{2}=\omega_{-k}$
and
$\omega_k+3k+1=\frac{k(3k-1)+6k+2}{2}=\frac{(k+1)(3k+2)}{2}=\omega_{k+1}$,
hence by replacing $k$ with $k+1$ we obtain
\begin{eqnarray*}
P&=&\sum_{n=-k}^{k+1} (-1)^n x^{\omega_n}\\
&&+(-1)^{k+1}x^{\omega_{k+1}+k+1}\Big(1+(1-x^{k+1})x^{k+1}+(1-x^{k+1})(1-x^{k+2})x^{2k+2}+\ldots\Big)
\end{eqnarray*}
which completes the induction, thus proving the pentagonal number theorem.

\section{Euler's recurrence relation for the divisor function}
\label{section:divisor}
We recall the divisor function $\int n$, the sum of all the divisors of $n$, including
itself.
Euler first discussed his recurrence relation for the divisor function
in
a letter to Goldbach 
on April 1, 1747,
Letter 113 in the Euler-Goldbach correspondence \cite{eulergoldbach},
OO827.
He begins the letter by saying, ``I have recently discovered
a very amazing order in the integers, which the sums of the divisors
of the natural numbers present, which appeared so much more peculiar
to me, since in this a great connection with the order
of the prime numbers appears to hide. Therefore I ask for some attention
to this.

If $n$ denotes any particular integral number, then $\int n$ should
denote the sum of all the divisors of this number $b$. Therefore we have:
\[
\begin{array}{l|l}
\int 1=1&\int 9=1+3+9=13\\
\int 2=1+2=3&\int 10=1+2+5+10=18\\
\int 3=1+3=4&\int 11=1+11=12\\
\int 4=1+2+4=7&\int 12=1+2+3+4+6+12=28\\
\int 5=1+5=6&\int 13=1+13=14\\
\int 6=1+2+3+6=12&\int 14=1+2+7+14=24\\
\int 7=1+7=8&\int 15=1+3+5+15=24\\
\int 8=1+2+4+8=15&\int 16=1+2+4+8+16=31\\
\multicolumn{2}{c}{\textrm{etc.}}
\end{array}
\]
Given this meaning for the symbol $\int$, I have found that
\begin{align*}
\int n=\int (n-1)+\int (n-2)-\int (n-5)-\int (n-7)+
\int (n-12)\\
+
\int (n-15)
-\int (n-22)-\int (n-26)+\textrm{etc.},
\end{align*}
where always the two signs $+$ and $-$ follow themselves. The order
of the derived numbers 1, 2, 5, 7, 12, 15, etc. 
appears from their differences, and if the same alternations are considered,
one immediately sees that,
{\footnotesize
\[
\begin{array}{p{0.25cm}p{0.25cm}p{0.25cm}p{0.25cm}p{0.25cm}p{0.25cm}p{0.25cm}p{0.25cm}p{0.25cm}p{0.25cm}p{0.25cm}p{0.25cm}p{0.25cm}p{0.25cm}p{0.25cm}p{0.25cm}p{0.25cm}p{0.25cm}p{0.25cm}p{0.25cm}p{0.25cm}}
&1,&&2,&&5,&&7,&&12,&&15,&&22,&&26,&&35,&&40,\\
\textrm{Diff.}&&1,&&3,&&2,&&5,&&3,&&7,&&4,&&9,&&5,&&11,\\
&51,&&57,&&70,&&77,&&92,&&100,&&117,&&126,&&145,&&\textrm{etc.}\\
\textrm{Diff.}&&6,&&13,&&7,&&15,&&8,&&17,&&9,&&19,&&\textrm{etc.}
\end{array}
\]
}
Furthermore it should be noted that in each case
one needs not take more terms once the negative numbers are come to,
and if such a term $\int 0$ appears, then for this the given number $n$
must be written, so that in such a case $\int 0=n$. The following
examples will illuminate the truth of this theorem:
\[
\begin{array}{rll}
&\textrm{When}&\textrm{then it will be}\\
1.&n=1;&\int 1=\int 0=1\\
2.&n=2;&\int 2=\int 1+\int 0=1+2=3\\
3.&n=3;&\int 3=\int 2+\int 1=3+1=4\\
4.&n=4;&\int 4=\int 3+\int 2=4+3=7\\
5.&n=5;&\int 5=\int 4+\int 3-\int 0=7+4-5=6\\
6.&n=6;&\int 6=\int 5+\int 4-\int 1=6+7-1=12\\
7.&n=7;&\int 7=\int 6+\int 5-\int 2-\int 0=12+6-3-7=8\\
8.&n=8;&\int 8=\int 7+\int 6-\int 3-\int 1=8+12-4-1=15\\
9.&n=9;&\int 9=\int 8+\int 7-\int 4-\int 2=15+8-7-3=13\\
10.&n=10;&\int 10=\int 9+\int 8-\int 5-\int 3=13+15-6-4=18\\
11.&n=11;&\int 11=\int 10+\int 9-\int 6-\int 4=18+13-12-7=12\\
12.&n=12;&\int 12=\int 11+\int 10-\int 7-\int 5+\int 0=12+18-8-6+12=28\\
&&\textrm{etc.}
\end{array}
\]
The reason for this order is 
not obvious, since one does not see how the numbers 1, 2, 5, 7, 12, 15, etc.
relate with the nature of the divisors. I can also not claim
that I have been able to give a rigorous proof of this either. However,
if I had no proof at all, one would still not be able to doubt the truth of it,
because over 300 cases always follow this rule. In the mean time, 
I have correctly derived this theorem from the following statement.

If $s=(1-x)(1-x^2)(1-x^3)(1-x^4)(1-x^5)$ etc., then also
$s=1-x-x^2+x^5+x^7-x^{12}-x^{15}+x^{22}+x^{26}-x^{35}-x^{40}+$ etc.,
where the exponents of $x$
are the same numbers which appeared earlier on;
and if this statement is true, which I do not doubt, despite the fact that I
do not have a rigorous demonstration, then the theorem is completely
justified. 

For from the double values of $s$, I obtain firstly
\[
\frac{ds}{s}=\frac{-dx}{1-x}-\frac{2xdx}{1-x^2}-\frac{3x^2dx}{1-x^3}-
\frac{4x^3dx}{1-x^4}-\frac{5x^4dx}{1-x^5}-\textrm{etc.}
\]
and then
\[
\frac{ds}{s}=\frac{-dx-2xdx+5x^4+7x^6-12x^{11}-15x^{14}+\textrm{etc.}}{1-x-xx+x^5+x^7-x^{12}-x^{15}+\textrm{etc.}}
\]
therefore we have
\begin{align*}
\frac{1+2x-5x^4-7x^6+12x^{11}+15x^{14}-
\textrm{etc.}}{1-x-x^2+x^5+x^7-x^{12}-x^{15}+\textrm{etc.}}=\\
\frac{1}{1-x}+\frac{2x}{1-x^2}+\frac{3x^2}{1-x^3}+\frac{5x^4}{1-x^5}+
\frac{4x^3}{1-x^4}+\frac{6x^5}{1-x^6}+\textrm{etc.}
\end{align*}
If however all of the last pieces are transformed into geometric
progressions, then one obtains for the same

{\footnotesize
\begin{tabular}{p{0.5cm}p{0.5cm}p{0.5cm}p{0.5cm}p{0.5cm}p{0.5cm}p{0.5cm}p{0.5cm}p{0.5cm}p{0.5cm}p{0.5cm}p{0.5cm}p{0.5cm}p{0.5cm}}
$1$&$+x$&$+x^2$&$+x^3$&$+x^4$&$+x^5$&$+x^6$&$+x^7$&$+x^8$&$+x^9$&$+x^{10}$&$+x^{11}$&$+x^{12}$&+etc.\\
&$+2x$&&$+2x^3$&&$+2x^5$&&$+2x^7$&&$+2x^9$&&$+2x^{11}$&&\\
&&$+3x^2$&&&$+3x^5$&&&$+3x^8$&&&$+3x^{11}$&&\\
&&&$+4x^3$&&&&$+4x^7$&&&&$+4x^{11}$&&\\
&&&&$+5x^4$&&&&&$+5x^9$&&&&\\
&&&&&$+6x^5$&&&&&&$+6x^{11}$&&\\
&&&&&&$+7x^6$&&&&&&&\\
&&&&&&&$+8x^7$&&&&&&\\
&&&&&&&&$+9x^8$&&&&&\\
&&&&&&&&&$+10x^9$&&&&\\
&&&&&&&&&&$+11x^{10}$&&&\\
&&&&&&&&&&&$+12x^{11}$&&\\
&&&&&&&&&&&&$+13x^{12}$&\\
&&&&&&&&&&&&&+etc.
\end{tabular}
}

that is:
{\footnotesize
\begin{align*}
1+\int 2x+\int 3x^2+\int 4x^3+\int 5x^4+\int 6x^5+\int 7x^6+\int 8x^7+\int 9x^8+
\int 10x^9+\textrm{etc.}\\
=
\frac{1+2x-5x^4-7x^6+12x^{11}+15x^{14}-22x^{21}-226x^{25}+35x^{34}+\textrm{etc.}}{1-x-x^2+x^5+x^7-x^{12}-x^{15}+x^{22}+x^{26}-x^{35}-\textrm{etc.}},
\end{align*}
}
from which the given theorem easily ensues. One sees however at the same time,
that this is not very obvious, and that without a doubt there are further
beautiful things hidden within this.'' 

The first paper in which Euler gives his recurrence relation for the divisor
function is
his ``D\'ecouverte d'une loi tout extraordinaire des nombres, par
rapport \`a la somme de leurs diviseurs'', presented to the Berlin Academy on June 22, 1747,
and published in the {\em Biblioth\`eque impartiale} in
1751 \cite{E175}, E175.
Euler proves this recurrence relation
in the same
way as he did in his April 1, 1747 letter to Goldbach, which we
explained above.

Euler writes next to d'Alembert about the pentagonal number theorem
in a letter on 
February 15, 1748, Letter 13 in the Euler-d'Alembert correspondence
in the {\em Opera Omnia} \cite{commercium5}, OO25.
Euler writes: ``Regarding this series that I spoke to you about, I 
found from it a very peculiar property about numbers with respect to the
sum of the divisors of each number. That $\int n$ represents
the sum of all the divisors of $n$ so that $\int 1=1; \int 2=3; \int 3=4; \int 4=7; \int 5=6; \int 6=12; \int 7=8$ etc. it seemed initially almost impossible
to discover any law in the sequence of numbers, but I found that each term
depends on some of the previous ones, according to this formula:
\[
\int n=\int (n-1)+\int (n-2)-\int (n-5)-\int (n-7)+\int (n-12)+\int (n-15)-
\int (n-22)-\textrm{etc.}
\]
where it is worthy of note $1^\circ$ that the numbers
\[
\begin{array}{cccccccccccccccccccc}
1,&&2,&&5,&&7,&&12,&&15,&&22,&&26,&&35,&&40,&\textrm{etc.}\\
&1&&3&&2&&5&&3&&7&&4&&9&&5&&
\end{array}
\]
are easily obtained by the differences considered alternately. 
$2^\circ$ In each case we only take
the numbers where the number after the $\int$ sign are non-negative.
$3^\circ$ If we obtain the term $\int 0$ or $\int (n-n)$ we will take $n$ as
the value.

Thus you will see that
\begin{align*}
\int 4=\int 3+\int 2=7;\\
\int 9=\int 8+\int 7-\int 4-\int 2=15+8-7-3=13;\\
\int 15=\int 14+\int 13-\int 10-\int 8+\int 3+\int 0=24+14-18-15+4+15=24;\\
\int 35=\int 34+\int 33-\int 30-\int 28+\int 23+\int 20-\int 13-\int 9+\int 0=\\
54+48-72-56+24+42-14-13+35=48.
\end{align*}
So every time that $n$ is a prime number we will find that $\int n=n+1$ and
since the nature of prime numbers enters this investigation, this law
seems to me even more remarkable.''

d'Alembert replies to Euler in a letter
on March 30, 1748, Letter 14 in the Euler-d'Alembert correspondence
in the {\em Opera omnia} \cite{commercium5}, OO26,
saying,
``That Sir, what comes to my mind while writing you, I only
have room left to say that your theorem on series seems very beautiful''.

On April 6, 1752, Euler presented his paper ``Observatio de summis
divisorum'' to the St. Petersburg Academy,
published in 1760 in the {\em Novi commentarii academiae scientiarum imperialis Petropolitanae} \cite{E243}, E243.
This paper just repeats what Euler said in the letter to Goldbach and paper
we considered above.

In the same
volume of this journal, Euler published the paper ``Demonstratio theorematis circa ordinem in summis
divisorum observatum'' \cite{E244}, E244, which again relates this recurrence
relation for the divisor function (although more clearly), and also
gives an inductive proof of the pentagonal number theorem, which we
discussed in Section~\ref{section:pentagonal}.

Goldbach writes a letter to Euler on May 9, 1752,
Letter 157 in the Euler-Goldbach correspondence \cite{eulergoldbach},
OO871,
in which he says that (Augustin Nathaniel) Grischow
has written to him about the presentation Euler gave to the
St. Petersburg Academy 
on the sums of divisors. Goldbach writes, ``I find myself at the present
time not in the position to judge, only you have insight into such matters,
though
let me have no doubt about the truth of everything that has been said in
this dissertation. In particular I saw with great pleasure 
that in the numbers 1, 2, 5, 7, 12, 15, 22, etc. was such a beautiful
order as was remarked before''.

\section{Euler's recurrence relation for the partition function}
We recall the partition function $n^{(m)}$, the number of ways
of expressing $n$ as a sum of positive integers less than or equal to $m$, disregarding order.
The only work in which Euler mentions 
his pentagonal number recurrence relation for the
partition function (there is no extant correspondence that discusses it)
is his paper ``De partitione numerorum'',
presented to
the St. Petersburg Academy on January 26, 1750 and published
in the {\em Novi commentarii academiae scientiarum imperialis Petropolitanae}
in 1753 \cite{E191}, E191. In \S 40, Euler says,
``Certainly
it is manifest from the nature of this matter for it truly to be a recurrent
series,  
with it 
arising from the expansion of this fraction:
\[
\frac{1}{(1-x)(1-x^2)(1-x^3)(1-x^4)(1-x^5)(1-x^6)\textrm{ etc.}}
\]
Therefore there will be a ladder
relation for this series, if the
denominator were expanded by multiplication. Indeed with this multiplication
having been carried out, the denominator will be found to be expressed
in the following way:
\[
1-x-x^2+x^5+x^7-x^{12}-x^{15}+x^{22}+x^{26}-x^{35}-x^{40}+x^{52}+x^{57}-x^{70}-x^{77}+\textrm{etc.}
\]
These powers of $x$ hold to such a rule,
from which formation is seen to be able to be determined with difficulty;
in the meantime however, from inspection it is soon apparent
for the pairs of terms to alternately be positive and negative.
No less the exponents of $x$ are observed to hold to a certain law, from
which the general term is gathered to be $x^\frac{n(3n \pm 1)}{2}$.
Namely no other powers occur aside from those whose exponents are
contained
in the form $\frac{3nn \pm n}{2}$, and 
for which the powers which arise from $n$ taken as an odd number
have the sign $-$, and indeed those formed from an even number, the sign $+$.''

Then in \S 41, Euler continues, ``This form therefore provides
to us the ladder of relation for the series which has been found,
which comes out to be
\begin{align*}
n^{(\infty)}=(n-1)^{(\infty)}+(n-2)^{(\infty)}-(n-5)^{(\infty)}-(n-7)^{(\infty)}+
(n-12)^{(\infty)}\\
+
(n-15)^{(\infty)}
-(n-22)^{(\infty)}-(n-26)^{(\infty)}
+
(n-35)^{(\infty)}+(n-40)^{(\infty)}\\-(n-51)^{(\infty)}
-(n-57)^{(\infty)}+
\textrm{etc.}
\end{align*}
Indeed by trying this rule of the progression
its place will easily be able to obtained.
For were it $n=30$ it will be found to be:
\[
30^{(\infty)}=29^{(\infty)}+28^{(\infty)}-25^{(\infty)}-
23^{(\infty)}+18^{(\infty)}+15^{(\infty)}-8^{(\infty)}-
4^{(\infty)}
\]
which indeed with these numbers taken from the table
\[
5604=4565+3718-1958-1255+385+176-22-5.
\]
And indeed in this way 
it will be pleasant for such series to always be continued''.

\section{Divergent series of the pentagonal numbers}
\label{section:series}
The last paper in which Euler discusses the pentagonal number theorem
and its applications
is his ``De mirabilibus proprietatibus numerorum pentagonalium'',
presented to the
St. Petersburg Academy on September 4, 1775,
and published in the {\em Acta academiae scientiarum imperialis Petropolitinae}
in 1783 \cite{E542}, E542.
In \S 2, Euler notes that every pentagonal number is one-third of a triangular
number (those numbers in the form $\frac{n(n+1)}{2}$), and in
\S 4 he gives his recurrence relation for the divisor function.
Then in \S 7,
Euler states the pentagonal number theorem, remarking that,
``This then deserves our admiration no less
than the properties mentioned above, with no fixed
rule apparent from which any connection can be
understood between the expansion of this
product and our pentagonal numbers''.

In the rest of this paper, Euler considers the summation
of series of the pentagonal numbers.
(We recall from Hans Rademacher's ``Comments on Euler's `De mirabilibus proprietatibus
numerorum pentagonalium' '' \cite{MR0262045}
that the Euler method of summation of a series
$\sum_{k=1}^\infty  a_k$ is defined by:
\[
\sum_{k=1}^\infty a_k=\lim_{N\to \infty} \sum_{m=1}^N \bigg(\frac{1}{2}\bigg)^{m+1}
\sum_{n=0}^m \binom{m}{n} a_n.
\] 
We find (cf. Rademacher \cite{MR0262045}) that by this summation, all the sums of the divergent series
which Euler gives in this paper may be rigorously justified.)

In \S 8, Euler states that
``Therefore with this series of powers of $x$ equal to this
infinite product, if it were set equal to nothing, so that
we have this equation:
\[
0=1-x^1-x^2+x^5+x^7-x^{12}-x^{15}+x^{22}+x^{26}-\textrm{etc.}
\]
it will involve all the roots, which the product equated to nothing
includes''.
Then in \S 9, Euler says, ``It is then clear for all the roots of each
power from unity to simultaneously be equal to the roots
of our equation'', and with $x^n=\cos{2i\pi} \pm \surd{-1}\sin{2i\pi}$,
``if for $n$ and $i$ are taken
all the successive integral numbers, the formula
$x=\cos{\frac{2i\pi}{n}} \pm \surd{-1} \frac{2i\pi}{n}$
will produce all the roots of our equation''. 

In \S 10, Euler states that for $\alpha,\beta,\gamma,\delta,\varepsilon$, etc.
the roots of unity, ``we gather for the sum of all these fractions
to be
$\frac{1}{\alpha}+\frac{1}{\beta}+\frac{1}{\gamma}+\frac{1}{\delta}+\textrm{etc.}=1$,
the sum of the products from two to be equal to $-1$, then indeed the sum
of the products from three to be equal to 0, the sum of the products
from four equal to 0, the sum of the products from five equal to -1,
the sum of the products from six equal to 0, the sum of the products from
seven equal to -1, etc. Then indeed we can also deduce the sum of all the squares of these fractions, namely
\[
\frac{1}{\alpha^2}+\frac{1}{\beta^2}+\frac{1}{\gamma^2}+\frac{1}{\delta^2}+\textrm{etc.}=3,
\]
the sum of the cubes
\[
\frac{1}{\alpha^3}+\frac{1}{\beta^3}+\frac{1}{\gamma^3}+\frac{1}{\delta^3}+
\textrm{etc.}=4,
\]
the sum of the biquadrates
\[
\frac{1}{\alpha^4}+\frac{1}{\beta^4}+\frac{1}{\gamma^4}+\frac{1}{\delta^4}+
\textrm{etc.}=7,
\]
and so on thusly, where however no order is obvious''. 

In \S \S 11-13, Euler notes a few other properties of
roots of unity: the reciprocal of a root
of unity is itself a root of unity;
an $n$th root
of unity $\alpha$ satisfies $\alpha^{in+\lambda}=\alpha^\lambda$ for $i$
an integer; and, roots of roots of unity are themselves roots of unity.

In \S 14, Euler writes that ``With us having assumed here $\alpha$
to be a root of the equation
$1-x^n=0$, we may run through the cases successively in which $n$ is 1, 2,
3, 4, etc. And foremost, if $n=1$ it is necessarily $\alpha=1$,
with which value having been substituted into our general equation
this form is induced:
\[
1-1-1+1+1-1-1+1+\textrm{etc.}
\]
which series is manifestly conflated from infinitely many periods, each of which
contains the terms $1-1-1+1$, from which the value of each period
is equal to 0, and thus the infinitely many periods
 which have been taken at once have a sum equal to 0''. 
In \S 15, Euler gives an alternate justification of this
with the sum of the Leibnitz (i.e. alternating) series
$1-1+1-1+1-1+\textrm{etc.}=\frac{1}{2}$. 

In \S 16, Euler writes ``We may now consider the case in which $n=2$ and
$\alpha \alpha=1$, in which indeed $\alpha$ is either $+1$ or $-1$.
We shall retain the letter $\alpha$ for designating either one of these
however, and with
\[
\alpha^3=\alpha, \, \alpha^4=1, \, \alpha^5=\alpha, \, \alpha^6=1, \, \textrm{etc.}
\]
having been substituted into our general equation this form will be induced:
\[
1-\alpha-1+\alpha+\alpha-1-\alpha-1|+1-\alpha-1+\alpha+\alpha-1-\alpha+1
\, \, \textrm{etc.}
\]
which series progresses evenly by certain periods, which are
replicated continuously, and each of which is composed of these eight
terms:
\[
1-\alpha-1+\alpha+\alpha-1-\alpha+1,
\]
of which the sum is 0, and thus is certain to vanish however large a number of integral
periods''. 
Then in \S 17, Euler states that the this series can be separated
into two subseries which must both be equal to 0:
\begin{align*}
1-1-1+1,+1-1-1+1,+1-1-1+1, \, \textrm{etc.}=0,\\
-\alpha+\alpha+\alpha-\alpha,-\alpha+\alpha+\alpha-\alpha,
-\alpha+\alpha+\alpha-\alpha, \, \textrm{etc.}=0.
\end{align*}

In \S \S 18 and 19, Euler considers the cases $\alpha^3=1$ and $\alpha^4=1$,
getting
\begin{align*}
-\alpha^2+\alpha^2+\alpha^2-\alpha^2,-\alpha^2+\alpha^2+\alpha^2-
\alpha^2,-\alpha^2+\alpha^2+\alpha^2-\alpha^2, \, \textrm{etc.}=0,\\
+\alpha^3-\alpha^3-\alpha^3+\alpha^3,+\alpha^3-\alpha^3-\alpha^3+\alpha^3,
+\alpha^3-\alpha^3-\alpha^3+\alpha^3, \, \textrm{etc.}=0.
\end{align*}
In \S 20, Euler examines $\alpha^5=1$, writing, ``here not all the powers
less then five occur. If it were $\alpha^5=1$, this periodic series
will be produced:
\[
\begin{array}{l|l}
1-\alpha+1-\alpha^2+\alpha^2-1+\alpha-1+\alpha^2-\alpha^2+1&-\alpha+1\\
-\alpha^2+\alpha^2-1+\alpha-1+\alpha^2-\alpha^2+1-\alpha&+1-\alpha^2+\alpha^2
\end{array}\, \textrm{etc.}
\]
where the powers $\alpha^3$ and $\alpha^5$ are entirely excluded''.
Then, in \S \S 20--22, Euler observes that for all higher roots
of unity the series are composed from periods which each sum to zero,
and in \S 23 that if $\alpha$ is an $n$th root of unity so
that $1-\frac{x}{\alpha}$ is a factor of $1-x^n$, then
$1-\frac{x}{\alpha}$ is also a factor of $1-x^{2n},1-x^{3n},
1-x^{4n}$, etc., and thus for each root of unity to have
infinite multiplicity in the equation
$\sum_{n=0}^\infty x^{\frac{3x^2 \pm x}{2}}=0$.

Euler considers divergent series of the pentagonal numbers
in \S \S 24--31 of the paper. In \S 24, 
he writes: ``We know 
moreover from the nature of equations, if an arbitrary
equation
\[
1+Ax+Bxx+Cx^3+Dx^4+\textrm{etc.}=0,
\]
should have two roots equal to $\alpha$, then in fact
for $\alpha$ to be a root of the equation born from
differentiating, namely:
\[
A+2Bx+3Cxx+4Dx^3+\textrm{etc.}=0,
\]
and if it has three roots equal to $\alpha$, then
in addition $\alpha$ will also be a root of the
equation born from differentiating after indeed we multiply
this differentiated equation by $x$
\[
1^2\cdot A+2^2\cdot Bx+3^2\cdot Cxx+4^2\cdot Dx^3+\textrm{etc.}=0,
\]
from which if this equation were to have $\lambda$ equal roots,
each of which were equal to $\alpha$, then it will always be
\[
1^\lambda \cdot A+2^\lambda \cdot B\alpha+3^\lambda \cdot C\alpha^3+
4^\lambda \cdot D\alpha^4+\textrm{etc.}=0,
\]
from which if for the grace of uniformity we were to multiply
this equation by $\alpha$, it will then be
\[
1^\lambda\cdot A\alpha+2^\lambda\cdot B\alpha^2+3^\lambda\cdot C\alpha^3+
4^\lambda\cdot D\alpha^4+\textrm{etc.}=0\textrm{''.}
\]

In \S 25, Euler writes: ``Therefore by putting $\alpha^n=1$
our equation formed from the pentagonal numbers
\[
1-x^1-x^2+x^5+x^7-x^{12}-x^{15}+\textrm{etc.}=0,
\]
shall have infinitely many roots equal to $\alpha$, and thus
$\alpha$ will be a root of all equations contained in this general
form:
\[
-1^\lambda x-2^\lambda x^2+5^\lambda x^5+7^\lambda x^7-12^\lambda
x^{12}-\textrm{etc.}=0
\]
for any integer whatsoever taken for $\lambda$. Therefore
it will always be
\[
-1^\lambda \alpha-2^\lambda \alpha^2+5\lambda \alpha^5+
7^\lambda \alpha^6-12^\lambda \alpha^{12}-\textrm{etc.}=0\textrm{''.}
\] 

In \S 26, Euler then writes, ``In order to make this clear,
we shall take $\alpha=1$, and it will always be
\[
-1^\lambda-2^\lambda+5^\lambda+7^\lambda-12^\lambda-15^\lambda+
\textrm{etc.}=0,
\]
and for the case $\lambda=0$ we have already probed the truth
of this equation. Were it therefore $\lambda=1$, it will be
revealed for the sum of this series diverging to infinity:
\[
-1-2+5+7-12-15+22+26-\textrm{etc.}
\]
to be equal to 0. Seeing moreover that this series is broken up,
that is, it is interpolated from two series, each of which may
be contemplated individually, by putting
\begin{align*}
s=-1+5-12+22-35+\textrm{etc. and}\\
t=-2+7-15+26-40+\textrm{etc.}
\end{align*}
where it ought to be shown to become $s+t=0$''.
In \S 27 Euler continues,
``Indeed from the doctrine
of series, which proceed with alternating signs, such as
$A-B+C-D+$ etc., it is known for the sum of this series
progressing into infinity to be equal to $\frac{1}{2}A-
\frac{1}{4}(B-A)+\frac{1}{8}(C-2B+A)-\frac{1}{16}(D-3C+3B-A)$ etc.,
which rule is thus conveniently related by differences,
namely by a rule based on the signs. From the series
of numbers $A,B,C,D,E$, etc. should be formed a series of differences,
so that each term of this series is subtracted from the following one:
it will be $a,b,c,d$, etc. Again, by the same law, from this series
of differences should be formed the series of second differences,
which shall be $a',b',c',d'$, etc.,
from this series in turn the series of third differences,
which shall be $a'',b'',c'',d'',e''$, etc.,
and thus in this way beyond, until constant differences
prevail. Then moreover from the first terms
of all these series the sum of the proposed series may thus
be determined, such that it will be
\[
\frac{1}{2}A-\frac{1}{4}a+\frac{1}{8}a'-\frac{1}{16}a''+
\frac{1}{32}a'''-\frac{1}{64}a''''+\textrm{etc.''}
\]

In \S 28, Euler then writes, ``With this rule which has been
established, with the signs having been switched it will 
be
\begin{align*}
-s=1-5+12-22+35-51+70-\textrm{etc. and}\\
-t=2-7+15-26+40-57+77-\textrm{etc.}
\end{align*}
These terms may be arranged in the following way and their differences
may be written:
{\footnotesize
\[
\begin{array}{p{0.1cm}p{0.1cm}p{0.1cm}p{0.1cm}p{0.1cm}p{0.1cm}p{0.1cm}
p{0.1cm}p{0.1cm}p{0.1cm}p{0.1cm}p{0.1cm}p{0.1cm}p{0.6cm}||p{0.1cm}
p{0.1cm}p{0.1cm}p{0.1cm}p{0.1cm}p{0.1cm}p{0.1cm}p{0.1cm}p{0.1cm}
p{0.1cm}p{0.1cm}p{0.1cm}p{0.1cm}p{2cm}}
1,&&5,&&12,&&22,&&35,&&51,&&70,&\textrm{ etc. }&2,&&7,&&15,&&26,
&&40,&&57,&&77,&\textrm{ etc.}\\
&4,&&7,&&10,&&13,&&16,&&19&&&&5,&&8,&&11,&&14,&&17,&&20&&\\
&&3,&&3,&&3,&&3,&&3,&&&&&&3,&&3,&&3,&&3,&&3,&&&\\
&&&0,&&0,&&0,&&0,&&&&&&&&0,&&0,&&0,&&0,&&&&
\end{array}
\]
}
Then it is therefore gathered to be
\begin{align*}
-s=\frac{1}{2}-\frac{4}{4}+\frac{3}{8}=-\frac{1}{8}, \, \textrm{that is} \,
s=\frac{1}{8}, \, \textrm{and in turn}\\
-t=\frac{2}{2}-\frac{5}{4}+\frac{3}{8}=\frac{1}{8}, \, \textrm{that is} \,
t=-\frac{1}{8}
\end{align*}
from which it is clearly
concluded to be $s+t=0$''.

In \S 29,
Euler writes,
``Although the rules by which these properties are supported
clearly leave no doubt, it will by no means be useless to exhibit further
the truth for the case $\lambda=2$ for it to be
\[
-1^2-2^2+5^2+7^2-12^2-15^2+22^2+\textrm{etc.}=0.
\]
Again this series is parted into two, which will be
with the sign changed:
\begin{align*}
s=+1^2-5^5+12^2-22^2+35^2-51^2+\textrm{etc.}\\
t=2^2-7^2+15^2-26^2+40^2-57^2+\textrm{etc.}
\end{align*}
and to find the sum of the first, the following operation
is instituted:
\[
\begin{array}{p{1.5cm}p{0.5cm}p{0.5cm}p{0.5cm}p{0.5cm}p{0.5cm}p{0.5cm}
p{0.5cm}p{0.5cm}p{0.5cm}p{0.5cm}p{0.5cm}p{0.5cm}p{0.5cm}}
\textrm{Series}&1,&&24,&&144,&&484,&&1225,&&2601,&&4900\\
\textrm{Diff. I}&&24,&&119,&&340,&&741,&&1376,&&2299&\\
\textrm{Diff. II}&&&95,&&221,&&401,&&635,&&923&&\\
\textrm{Diff. III}&&&&126,&&180,&&234&&288&&&\\
\textrm{Diff. IV}&&&&&54,&&54,&&54&&&&\\
\textrm{Diff. V}&&&&&&0,&&0&&&&&
\end{array}
\]
Then it will therefore be
\[
s=\frac{1}{2}-\frac{24}{4}+\frac{95}{8}-\frac{126}{16}+\frac{54}{32}=+\frac{3}{16}.
\]
In a similar way for the other series
\[
\begin{array}{p{1.5cm}p{0.5cm}p{0.5cm}p{0.5cm}p{0.5cm}p{0.5cm}p{0.5cm}p{0.5cm}p{0.5cm}p{0.5cm}p{0.5cm}p{0.5cm}p{0.5cm}p{0.5cm}}
\textrm{Series}&4,&&49,&&225,&&676,&&1600,&&3249,&&5929\\
\textrm{Diff. I}&&45,&&176,&&451,&&924,&&1649,&&2680&\\
\textrm{Diff. II}&&&131,&&275,&&473,&&725,&&1031&&\\
\textrm{Diff. III}&&&&144,&&198,&&252,&&306&&&\\
\textrm{Diff. IV}&&&&&54,&&54,&&54&&&&\\
\textrm{Diff. V}&&&&&&0,&&0&&&&&
\end{array}
\]
Then it may be concluded
\[
t=\frac{4}{2}-\frac{45}{4}+\frac{131}{8}-\frac{144}{16}+\frac{54}{32}=
-\frac{3}{16}.
\]
From this it prevails for the total sum to become $s+t=0$.''

In \S 30, ``We shall now consider as well the square roots,
that is when it is $\alpha^2=1$, and then such a series will arise:
\[
-1^\lambda \cdot \alpha-2^\lambda+5^\lambda\cdot \alpha-7^\lambda\cdot \alpha-
12^\lambda-15^\lambda \cdot \alpha+22^\lambda+26^\lambda-\textrm{etc.}=0,
\]
from which if we separate the terms containing unity and $\alpha$ 
from each other, we shall obtain two series equal to nothing, namely:
\[
-2^\lambda-12^\lambda+22^\lambda+26^\lambda-40^\lambda-70^\lambda+92^\lambda+
\textrm{etc.}=0
\]
and
\[
-1^\lambda\cdot \alpha+5^\lambda \cdot \alpha+7^\lambda\cdot \alpha-15^\lambda\cdot \alpha-
35^\lambda\cdot \alpha+51^\lambda\cdot \alpha+57^\lambda\cdot \alpha-
\textrm{etc.}=0.
\]
If indeed we want to display the truth of these series in the same way
which we did before, each ought to be divided
into four others, until in the end we reach constant differences. 
And indeed, if this work were undertaken, it will be able to be
certain, for the aggregate of all the parts to be equal to 0''.

In \S 31, Euler says,
``Now most generally the total problem
is embraced, and it may be $\alpha^n=1$, and we shall search for the
series which contains all the powers $\alpha^r$.
To this end, from all our pentagonal numbers we shall pick out
those which divided by $n$ leave the very residue $r$. Therefore
were these pentagonal numbers $A,B,C,D,E$, etc., namely all those
of the form $\gamma n+r$, and the sign of which $\pm$, which will agree with
these, may be noted with care. Then indeed it will always be
\[
\pm A^\lambda \pm B^\lambda \pm C^\lambda \pm D^\lambda \pm
\textrm{etc.}=0,
\]
where any integral value may be taken for the exponent $\lambda$.
Then all the series which we have elicited so far, and of which we have
shown for the sums to be equal to nothing, are contained
in this most general form.''

\section{Conclusions}
\label{section:conclusion}
As we noted in in Section~\ref{section:introduction},
the pentagonal number theorem is a special case of the Jacobi triple
product identity. Using
the Jacobi triple product identity
we can also obtain
an identity for the cube of the product in \eqref{equation:maintheorem},
$\prod_{m=1}^\infty (1-x^m)^3=\sum_{n=0}^\infty
(-1)^n(2n+1)x^{\frac{n(n+1)}{2}}$. 
C. G. J. Jacobi observes this in his article ``Note sur les
fonctions elliptiques'' in 1828 \cite{jacobi}, saying that he
always found Euler's
pentagonal number theorem ``a very surprising and admirable feat''.
Indeed, the pentagonal number theorem is one of the first results in
the theory of $q$-series and theta functions.

Euler considered the series
expansion of more general infinite products
such as 
$(1+mz)(1+m^2z)(1+m^3z)(1+m^4z)$ etc. (for example in
his November 10, 1742 letter to Niklaus I Bernoulli, his
April 5, 1746 letter to Goldbach, and \S 18 of ``Observationes analyticae variae
de combinationibus'' \cite{E158}, which we discussed in
Section~\ref{section:background}),
but he does not seem to have found the Jacobi triple product identity.

In terms of the product representation of theta functions, the pentagonal number
theorem is
$\vartheta_4(q/6;q)=(q;q)_\infty$ \cite[\S 78, Chapter 10]{MR0364103}, 
for $q^{2/3}=x$.

It would be good to further explore Euler's work on partitions.
Kiselev and Matvievskaja
discuss \cite{MR0199077} unpublished notes of Euler on partitions.

\section{Acknowledgements}
We would like to thank Dr. Paul Mezo for translating parts of
Euler's correspondence from the German, Christian L\'eger for
translating parts of Euler's correspondence from the French,
and Dr. Kenneth S. Williams for helpful discussions about Euler's
divergent series of the pentagonal numbers.

\bibliographystyle{apalike}
\bibliography{pentagonal}

\end{document}